\documentclass[aps,pre,reprint,superscriptaddress,longbibliography,amsmath,amssymb]{revtex4-2}
\usepackage[T1]{fontenc}
\usepackage[utf8]{inputenc}
\usepackage{amsthm}
\usepackage[hidelinks]{hyperref}
\hypersetup{
  pdftitle={Local Clustering Tests and Botnet Localization in Sparse Random Geometric Networks},
  pdfauthor={Luca Di Persio, Gianmarco Bet, Youness Outaleb},
  pdfkeywords={random geometric graph, botnet detection, local clustering coefficient, triangle statistic, scan test, exact localization.}
}

\theoremstyle{plain}
\newtheorem{theorem}{Theorem}[section]
\newtheorem{lemma}[theorem]{Lemma}
\newtheorem{proposition}[theorem]{Proposition}
\newtheorem{corollary}[theorem]{Corollary}
\theoremstyle{definition}

\theoremstyle{remark}

\newcommand{\Bin}{\mathrm{Bin}}
\newcommand{\Vol}{\mathrm{Vol}}
\newcommand{\Var}{\mathrm{Var}}
\newcommand{\Cov}{\mathrm{Cov}}
\newcommand{\TT}{\mathbb{T}}
\newcommand{\RR}{\mathbb{R}}
\newcommand{\ZZ}{\mathbb{Z}}
\newcommand{\Prob}{\mathbb{P}}
\newcommand{\Exp}{\mathbb{E}}
\newcommand{\one}{\mathbf{1}}
\newcommand{\Cd}{C_d}
\newcommand{\thetad}{\theta_d}
\newcommand{\mun}{\mu_n}
\newcommand{\distT}{D_{\TT}}
\newcommand{\Bh}{\widehat B}
\newcommand{\Risk}{\mathcal{R}}
\newcommand{\calB}{\mathcal{B}}

\begin{document}

\title{Local Clustering Tests and Botnet Localization in Sparse Random Geometric Networks}
\author{Luca Di Persio}
\affiliation{Department of Computer Science, University of Verona, Verona, Italy}
\author{Gianmarco Bet}
\affiliation{Department of Mathematics and Computer Science ``Ulisse Dini'', University of Florence, Florence, Italy}
\author{Youness Outaleb}
\affiliation{Department of Mathematics, University of Trento, Trento, Italy}

\begin{abstract}
We use local clustering to test for botnets in sparse random geometric networks. Under the null hypothesis, the vertices are independent uniform points on a fixed-dimensional torus, with edges determined by a distance threshold. Under the alternative, an unknown set of vertices ignores the geometry and connects independently with the same marginal edge probability. Every prescribed vertex has the same marginal binomial degree law under both hypotheses, so its degree alone carries no information. Triangle closure, however, distinguishes the two models. Under the null, conditional on a degree of at least two, the local clustering coefficient has exact mean $\Cd>0$, a known constant that depends only on the dimension. Under the alternative, the corresponding mean is $p_n$ for a bot. For a non-bot, the corresponding mixture mean converges to $\Cd$ when the botnet size is $o(n)$. We prove non-asymptotic concentration bounds, an exponential bound on the worst-case risk of a scan test, and exact-support recovery. For every deterministic sequence $1\leq k_n=o(n)$, if $\mun=(n-1)p_n$ satisfies $\mun/\log n\to\infty$, the local rule recovers the botnet with probability tending to one, uniformly over all label sets of size $k_n$. The proof combines geometric U-statistics, bounded differences on random neighborhoods, and a stability bound showing that genuine vertices retain their geometric clustering when bots form only a vanishing fraction of their neighbors. The test does not use the latent positions. The theorem assumes that the model parameters are known, although $\Cd$ can instead be estimated consistently. In the high-expected-degree regime, the rule also works for a single bot.
\end{abstract}

\maketitle

\section{Introduction}
Random geometric graphs describe networks whose edges reflect latent distance. Gilbert introduced a random plane network by placing points in space and joining nearby pairs \cite{gilbert1961random}. Penrose developed the probabilistic theory of their components, connectivity, subgraph counts, degrees, and limit laws \cite{penrose2003random}. Latent space models use a softer version of the same idea: vertices have unobserved positions, and their connection probabilities depend on latent distance \cite{hoff2002latent}. After the latent positions are integrated out, these vertex-exchangeable models generally have dependent edge indicators. Conditional independence given the positions is not contradicted by this marginal dependence.

We use local transitivity. Two neighbors of one vertex have a positive probability of being close to each other, and the clustering coefficient of Watts and Strogatz measures this closure \cite{watts1998collective}. For geometric graphs, Dall and Christensen computed the coefficient and showed that it differs from the Erd\H{o}s--R\'enyi value \cite{dall2002random}. Dhara, van Leeuwaarden, and Mukherjee later expressed it in terms of the incomplete beta function \cite{dhara2018corrected}. In fixed dimension, the coefficient does not vanish as the connection radius shrinks, so it remains a local measure of the geometry.

Bet, Bogerd, Castro, and van der Hofstad introduced botnet detection in random geometric networks \cite{bet2020detecting}. They consider a single observed graph, which is random geometric under the null. Under the alternative, a small unknown set of bot vertices connects independently with the same marginal edge probability as a geometric edge. Consequently, every individual edge indicator has the same marginal law under the two hypotheses, as does the degree of any prescribed vertex. The higher-order dependence among the edge indicators nevertheless differs. The tests in \cite{bet2020detecting} are based on isolated stars and average distances. Their isolated-star procedure also localizes bots in a setting where $k_n$ is known. Under the stated regimes, its expected normalized symmetric-difference loss tends to zero. This gives $o(k_n)$ expected label errors when $k_n\to\infty$ and exact-set recovery when $k_n$ is fixed. Here we instead prove exact-support recovery, $P_B(\widehat B_n=B)\to1$, for a clustering rule under a stronger high-degree assumption. The isolated-star test exploits local patterns created by random bot attachments, whereas the distance test exploits the shorter graph distances produced by non-geometric long-range links.

We use the local clustering coefficient of each vertex: the fraction of pairs of its neighbors that are adjacent. Under the null, the neighbors of a vertex lie in a small metric ball, and its degree-conditional local clustering coefficient has exact mean $\Cd$. A bot's neighbors do not share a common latent location; conditional on the bot's degree, two of them are adjacent with probability $p_n$. A genuine vertex under the alternative has a mixture center because some of its neighbors may be bots. This center converges to $\Cd$ when $k_n=o(n)$. Since $p_n\to0$, the gap remains of constant order and is large enough to support a scan over all vertices when the expected degree $\mun$ grows faster than $\log n$.

Random geometric graph theory supplies both the null model and tools for handling spatial dependence \cite{penrose2003random}. Clustering is a standard local network descriptor \cite{watts1998collective}, and its connection with latent geometry has also been studied in other network models \cite{krioukov2016clustering}. Signed-triangle statistics have been used to test for latent geometry \cite{bubeck2016testing}, while local clustering has been used as an empirical feature for botnet detection \cite{chowdhury2017botnet}. Concentration for geometric subgraph-count U-statistics has been studied in Poisson models \cite{bachmann2018concentration}. Bet et al. formulate detection against a geometry-breaking alternative as a minimax testing problem \cite{bet2020detecting}. More broadly, graph anomaly detection studies tests based on local subgraph features \cite{akoglu2015graph}. To our knowledge, this is the first exact-support result for the Bet et al. model based on a thresholded vertex local clustering coefficient; the guarantee is uniform over all label sets of size $k_n=o(n)$.

Under the null, we derive an exact scale-free conditional law and a conditional variance formula based on $\Cd$. Under the alternative, we compute the exact finite-$n$ conditional center of a genuine vertex. Two-sided concentration bounds for null vertices and bots then give explicit exponential bounds on the worst-case risk of the scan. For every deterministic sequence $1\leq k_n=o(n)$ satisfying $\mun/\log n\to\infty$, the vertices with typical degree and abnormally low local clustering are exactly the bots with probability tending to one, uniformly over their labels. This is exact-support recovery, not just detection. The proof uses the fact that bots form a vanishing fraction of a genuine vertex's neighbors.

\section{Main results}
Let $d\geq 2$ be fixed. The vertex set is $V_n=\{1,\ldots,n\}$. The latent space is the unit torus $\TT^d=[0,1]^d$ equipped with the torus metric
$$
\distT(x,y)=\min_{z\in\ZZ^d}\|x-y+z\|_2.
$$
Let $X_1,\ldots,X_n$ be independent uniform points on $\TT^d$. Let $p_n\in(0,1)$, let $p_n\to0$, and define
$$
\mun=(n-1)p_n.
$$
The radius $r_n$ is chosen by
$$
\Prob\{\distT(X_1,X_2)\leq r_n\}=p_n.
$$
Since $p_n\to0$, for all sufficiently large $n$ one has $r_n < \frac{1}{3}$.

For such $n$, the torus ball of radius $r_n$ lifts to a Euclidean ball and
$$
p_n=\thetad r_n^d,
\qquad
\thetad=\Vol\{x\in\RR^d:\|x\|_2\leq1\}.
$$
Unless noted otherwise, we assume that $d$ and $p_n$ are known; for the tests, this means that $\Cd$ and $\mun$ are known. The latent positions are not observed. Whenever we condition on $X_i=x$, we use the translation-invariant regular conditional distribution induced by the product density.

Under the null hypothesis $H_0$, the observed graph $G_n=(V_n,E_n)$ is the random geometric graph
$$
\{i,j\}\in E_n
\quad\Longleftrightarrow\quad
\distT(X_i,X_j)\leq r_n.
$$
Write $P_0$ for its law. Under the alternative, an unknown set $B\subset V_n$ with $|B|=k_n$, where $1\leq k_n\leq n$, is the botnet. Edges between two non-bot vertices are still generated by the geometric rule. Every edge incident to at least one bot is generated independently with probability $p_n$, independently of all locations and of all other bot-incident edges. Write $P_B$ for the law under a fixed botnet $B$. The class of alternatives is
$$
\calB_{n,k_n}=\{B\subset V_n: |B|=k_n\}.
$$
A test is a measurable function $\psi_n(G_n)\in\{0,1\}$, with value one indicating rejection of $H_0$. Its worst-case risk is
$$
\Risk_n(\psi_n)
=P_0(\psi_n=1)+\max_{B\in\calB_{n,k_n}}P_B(\psi_n=0).
$$
A sequence of tests is asymptotically powerful when $\Risk_n(\psi_n) \to 0$.

For a vertex $i$, write
$$
N(i)=\{j\in V_n\setminus\{i\}:\{i,j\}\in E_n\},
\qquad
D_i=|N(i)|.
$$
The number of triangles incident to $i$ is
$$
T_i=
\sum_{\{j,k\}\subset V_n\setminus\{i\}}
\one\{\{i,j\}\in E_n,\{i,k\}\in E_n,\{j,k\}\in E_n\}.
$$
The local clustering coefficient is
$$
C_i=
\begin{cases}
\frac{T_i}{\binom{D_i}{2}},& D_i\geq2,\\
0,& D_i<2.
\end{cases}
$$
For $r_n<1/3$ and distinct vertices $i,j,k$, define the geometric clustering constant \cite{bet2020detecting} by
$$
\Cd
=P_0\bigl(\{j,k\}\in E_n\mid \{i,j\}\in E_n,\ \{i,k\}\in E_n\bigr).
$$
Dall and Christensen give the equivalent overlap formula \cite{dall2002random}:
$$
\Cd=\frac{1}{\thetad^2}
\int_{B(0,1)}\Vol\bigl(B(0,1)\cap B(x,1)\bigr)\,dx.
$$
Here $\thetad=\Vol(B(0,1))$, and $B(x,r)$ is the closed Euclidean ball in $\RR^d$ with center $x$ and radius $r$. Define
$$
I_z(a,b)=\frac{\int_0^z t^{a-1}(1-t)^{b-1}\,dt}
{\int_0^1 t^{a-1}(1-t)^{b-1}\,dt}
$$
which is the regularized incomplete beta function. With $a=(d+1)/2$, the formulas in \cite{dhara2018corrected,bet2020detecting} give
$$
\begin{aligned}
&\Cd =d\int_0^1 s^{d-1}I_{1-s^2/4}\left(a,\frac12\right)\,ds\\
&\quad =I_{3/4}\left(a,\frac12\right) +I_{1/4}(a,a).
\end{aligned}
$$
The number $\Cd$ lies strictly between zero and one.

\begin{theorem}\label{thm:constant}
For all sufficiently large $n$, every vertex $i$, every set $S\subset V_n\setminus\{i\}$ with $|S|=m\geq2$, and every location $x\in\TT^d$, one has
$$
\Exp_0\left[C_i\mid X_i=x,\, N(i)=S\right]=\Cd.
$$
Consequently,
$$
\Exp_0\left[C_i\mid D_i=m\right]=\Cd,
\; \; 2\leq m\leq n-1.
$$
\end{theorem}

\begin{proof}
Fix $i$, $S$, $m$, and $x$. Given $X_i = x$ and $N(i) = S$, the points $(X_j)_{j \in S}$ are independent and uniform on $B_{\TT}(x,r_n)$. Indeed, the conditioning places the points indexed by $S$ inside the ball and the other points outside it. Write the neighbor locations as $Y_1,\ldots,Y_m$. Then
$$
C_i=\frac{1}{\binom{m}{2}} \sum_{1\leq a<b\leq m}\one\{\distT(Y_a,Y_b) \leq r_n\}.
$$
Taking expectation gives
$$
\begin{aligned}
&\Exp_0[C_i\mid X_i=x,N(i)=S]\\
&\quad =\frac{1}{\binom{m}{2}}\sum_{1\leq a<b\leq m}
\Prob_0\{\distT(Y_a,Y_b)\leq r_n\}.
\end{aligned}
$$
All summands have the same expectation, so
$$
\Exp_0\left[C_i\mid X_i=x,\,N(i)=S\right]
=\Prob\{\distT(Y_1,Y_2)\leq r_n\}.
$$
Lemma \ref{app:torus-lifting} lets us translate the torus ball to the origin and lift it to $B(0,r_n)$. Uniformity and Fubini's theorem then give
$$
\begin{aligned}
&\Prob\{\distT(Y_1,Y_2)\leq r_n\}\\
&\quad =\frac{1}{\thetad^2r_n^{2d}}
\int_{B(0,r_n)}
\Vol\bigl(B(0,r_n)\cap B(u,r_n)\bigr)\,du\\
&\quad =\frac{1}{\thetad^2}
\int_{B(0,1)}
\Vol\bigl(B(0,1)\cap B(v,1)\bigr)\,dv =\Cd,
\end{aligned}
$$
The second equality uses $u=r_nv$. Since the result is independent of $x$ and $S$, conditioning only on $D_i=m$ gives the last assertion.
\end{proof}

When $r_n<1/3$, the whole conditional law, not only its mean, is independent of the connection scale.

\begin{proposition}\label{prop:scale-free}
Let $U_1,\ldots,U_m$ be independent uniform random variables on $B(0,1)$, and define
$$
\mathcal C_{m,d} = \frac{1}{\binom{m}{2}} \sum_{1\leq a<b\leq m}\one\{\|U_a-U_b\|_2\leq1\}.
$$
For all sufficiently large $n$, every vertex $i$, and every $2\leq m\leq n-1$, the conditional law of $C_i$ given $D_i=m$ under $P_0$ is the law of $\mathcal C_{m,d}$.
\end{proposition}

\begin{proof}
Given $X_i = x$ and a neighbor set $N(i) = S$ of size $m$, write the neighbor locations as $Y_1,\ldots,Y_m$. As in the proof of Theorem \ref{thm:constant}, they are independent and uniform on $B_{\TT}(x,r_n)$. Translate by $-x$ and divide by $r_n$. For $r_n < \frac{1}{3}$, this gives independent uniform points $U_1,\ldots,U_m$ in $B(0,1)$, and the torus-lifting argument gives, for every pair,
$$
\distT(Y_a,Y_b)\leq r_n
\quad\Longleftrightarrow\quad
\|U_a-U_b\|_2\leq1
$$
The displayed equivalence shows that the conditional statistic has law $\mathcal C_{m,d}$ for every $x$ and $S$, independently of $n$ and $p_n$. Averaging over $x$ and $S$ conditional on $D_i=m$ proves the proposition.
\end{proof}

\begin{proposition}\label{prop:variance}
Let $U,V,W$ be independent uniform random variables on $B(0,1)$. Define
$$
h(u,v)=\one\{\|u-v\|_2\leq1\},
\qquad
q(u)=\Prob\{\|u-V\|_2\leq1\}.
$$
Set
$$
\sigma_{1,d}^2=\Var(q(U)),
\qquad
\sigma_{2,d}^2=\Cd(1-\Cd).
$$
For every vertex $i$, all sufficiently large $n$, and every feasible $2\leq m\leq n-1$,
$$
\Var_0\left(C_i\mid D_i=m\right)
=\frac{4(m-2)}{m(m-1)}\sigma_{1,d}^2
+\frac{2}{m(m-1)}\sigma_{2,d}^2.
$$
In particular
$$
\Var_0\left(C_i\mid D_i=m\right)
\leq \frac{4\sigma_{1,d}^2+2\sigma_{2,d}^2}{m-1}.
$$
\end{proposition}

\begin{proof}
This is the covariance calculation for an order-two U-statistic \cite{hoeffding1948class}. After translation and scaling, the conditional neighbor locations are independent uniform points $U_1,\ldots,U_m$ in $B(0,1)$, so
$$
C_i=\binom{m}{2}^{-1}\sum_{1\leq a<b\leq m}h(U_a,U_b).
$$
Let $N_m=\binom{m}{2}$. Expanding the square gives
$$
\begin{aligned}
&\Var(C_i\mid D_i=m)\\
&\quad =N_m^{-2}\sum_{a<b}\sum_{c<e}\Cov(h(U_a,U_b),h(U_c,U_e)).
\end{aligned}
$$
Disjoint index-pairs give zero covariance. Coincident pairs contribute $\Var(h(U,V))=\Cd(1-\Cd)=\sigma_{2,d}^2$. Pairs sharing one index contribute
$$
\Cov(h(U,V),h(U,W))
=\Exp[q(U)^2]-\Cd^2
=\sigma_{1,d}^2.
$$
For each of the $N_m$ index-pairs, there are $2(m-2)$ ordered choices of a second pair sharing one endpoint. Therefore
$$
\Var(C_i\mid D_i=m)
=N_m^{-2}\{N_m\sigma_{2,d}^2+N_m2(m-2)\sigma_{1,d}^2\},
$$
which is the asserted identity after substituting $N_m=m(m-1)/2$.
\end{proof}

The variance formula alone does not control a scan over all vertices. We use the following bounded-differences estimate \cite{mcdiarmid1989method,boucheron2013concentration}.

\begin{lemma}\label{lem:null-conc}
For every vertex $i$, all sufficiently large $n$, every integer $2\leq m\leq n-1$, and every $\varepsilon>0$,
$$
P_0\left(|C_i-\Cd|\geq\varepsilon\mid D_i=m\right)
\leq2\exp\left(-\frac{\varepsilon^2m}{2}\right).
$$
\end{lemma}

\begin{proof}
Condition on $X_i=x$ and $N(i)=S$, where $|S|=m$. Let $Y_1,\ldots,Y_m$ be the conditional neighbor locations. Conditional on this event, the local clustering coefficient can be written as the function
$$
\begin{aligned}
&C_i = f(Y_1,\ldots,Y_m)\\
&\quad =\binom{m}{2}^{-1}\sum_{1\leq a<b\leq m}\one\{\distT(Y_a,Y_b)\leq r_n\}.
\end{aligned}
$$
Changing one coordinate $Y_\ell$ can change only the $m-1$ indicators involving that coordinate. Since the normalizing factor is $\binom{m}{2}^{-1}$, the absolute change in $f$ is at most
$$
\frac{m-1}{\binom{m}{2}}=\frac{2}{m}.
$$
Since $Y_1,\ldots,Y_m$ are independent under this conditional law, McDiarmid's inequality, stated as Lemma \ref{app:mcdiarmid}, gives
$$
\begin{aligned}
&P_0\left(|f-\Exp[f\mid X_i=x,\!N(i)=S]|\geq\varepsilon\mid X_i=x,\!N(i)=S\right)\\
&\quad \leq2\exp\left(-\frac{2\varepsilon^2}{m(2/m)^2}\right).
\end{aligned}
$$
The conditional mean is $\Cd$ by Theorem \ref{thm:constant}. Averaging over $x$ and $S$ conditional on $D_i=m$ proves the lemma.
\end{proof}

For a bot, the conditional center is $p_n$. The following bound is uniform in the choice of botnet.

\begin{lemma}\label{lem:bot-conc}
Fix a botnet $B\in\calB_{n,k_n}$ and a bot vertex $b\in B$. For every integer $2\leq m\leq n-1$, we have
$$
\Exp_B\left[C_b\mid D_b=m\right]=p_n.
$$
For $m<2$, the conditional mean is zero. Moreover, under $P_B$, for every integer $2\leq m\leq n-1$ and every $\varepsilon>0$,
$$
P_B\left(|C_b-p_n|\geq\varepsilon\mid D_b=m\right)
\leq2\exp\left(-\frac{\varepsilon^2m}{4}\right).
$$
\end{lemma}

\begin{proof}
Fix $S\subset V_n\setminus\{b\}$ with $|S|=m$ and condition on $N(b)=S$. Let $S_0=S\setminus B$ and $S_1=S\cap B$. Because the edges incident to $b$ are independent of all locations, the vertices in $S_0$ remain independent and uniform on $\TT^d$. Edges within $S$ that touch $S_1$ are independent Bernoulli variables with parameter $p_n$; they are also independent of these locations and of the event $N(b)=S$. Edges with both endpoints in $S_0$ follow the geometric rule.

Every pair $\{u,v\}\subset S$ is therefore joined with conditional probability $p_n$. This is true by construction if the pair meets $S_1$; if both vertices lie in $S_0$, it follows from their independent uniform locations and the definition of $r_n$. Hence
$$
\Exp_B\left[C_b\mid N(b)=S\right]
=\binom{m}{2}^{-1}\sum_{\{u,v\}\subset S}p_n
=p_n.
$$
Averaging over $S$ conditional on $D_b = m$ gives the stated mean.

For concentration, regard $C_b$ as a function of the locations in $S_0$ and of the Bernoulli variables for pairs in $S$ that meet $S_1$. Changing one genuine location can change at most $m - 1$ edges inside the neighbor set and therefore changes $C_b$ by at most $2/m$. Changing one Bernoulli edge variable changes $C_b$ by $\binom{m}{2}^{-1}$. The sum of squared bounded-difference constants is at most
$$
|S_0|\left(\frac{2}{m}\right)^2
+\binom{m}{2}\left(\binom{m}{2}^{-1}\right)^2
\leq \frac{4}{m}+\frac{2}{m(m-1)}
\leq \frac{8}{m}.
$$
McDiarmid's inequality gives
$$
\begin{aligned}
&P_B\left(|C_b-p_n|\geq\varepsilon\mid N(b)=S\right) \leq2\exp\left(-\frac{2\varepsilon^2}{8/m}\right)\\
&\quad =2\exp\left(-\frac{\varepsilon^2m}{4}\right).
\end{aligned}
$$
Averaging over $S$ proves the bound.
\end{proof}

\paragraph{Mean at a genuine vertex under the alternative.}
Bot neighbors shift the mean of a genuine vertex away from $\Cd$. For all sufficiently large $n$ (so that $r_n<1/3$), if $i\notin B$ and $2\leq m\leq n-1$, then
$$
\Exp_B[C_i\mid D_i=m]
=p_n+(\Cd-p_n)
\frac{\binom{n-k_n-1}{2}}{\binom{n-1}{2}}.
$$
Given $X_i=x$, the $n-1$ incident edge indicators are independent Bernoulli variables with parameter $p_n$. Thus, conditional on $D_i=m$, the set $N(i)$ is a uniform $m$-subset of $V_n\setminus\{i\}$. A sampled pair of neighbors closes with probability $\Cd$ when both vertices are genuine and with probability $p_n$ when at least one is a bot. The factor in the display is exactly the probability that two labels sampled without replacement from $V_n\setminus\{i\}$ are genuine. Therefore
$$
\begin{aligned}
&0 \leq \Cd-\Exp_B[C_i\mid D_i=m] =(\Cd-p_n)\frac{k_n(2n-3-k_n)}{(n-1)(n-2)}\\
&\quad =O\left(\frac{k_n}{n}\right)
\end{aligned}
$$
whenever $p_n<\Cd$ and $k_n=o(n)$. The localization proof instead uses $C_i^g$ and the deterministic contamination bound below.

\begin{theorem}\label{thm:nonasymp}
Fix $\eta\in(0,1)$ and $\tau\in(0,\Cd)$, and put $\mun=(n-1)p_n$. Define the test
$$
\psi_{n,\eta,\tau}
=
\one\{\exists i\in V_n:D_i\geq(1-\eta)\mun\text{ and }C_i\leq\tau\}.
$$
Assume $r_n<\frac{1}{3}$, $p_n<\tau<\Cd$, and $(1-\eta)\mun\geq2$.
$$
P_0(\psi_{n,\eta,\tau}=1)
\leq
2n\exp\left(-\frac{(\Cd-\tau)^2(1-\eta)\mun}{2}\right).
$$
For every integer botnet size $1\leq k_n\leq n$,
$$
\begin{aligned}
&\max_{B\in\calB_{n,k_n}}P_B(\psi_{n,\eta,\tau}=0) \leq
\exp\left(-\frac{\eta^2\mun}{2}\right)\\
&\quad +2\exp\left(-\frac{(\tau-p_n)^2(1-\eta)\mun}{4}\right).
\end{aligned}
$$
Consequently
$$
\begin{aligned}
&\Risk_n(\psi_{n,\eta,\tau}) \leq 2n\exp\left(-\frac{(\Cd-\tau)^2(1-\eta)\mun}{2} \right)
\quad\\
&\quad +\exp\left(-\frac{\eta^2\mun}{2}\right) +2\exp\left(-\frac{(\tau-p_n)^2(1-\eta)\mun}{4} \right).
\end{aligned}
$$
\end{theorem}

\begin{proof}
Let
$$
E_i=\{D_i\geq(1-\eta)\mun,\ C_i\leq\tau\}.
$$
Under $H_0$, exchangeability and the union bound yield
$$
P_0(\psi_{n,\eta,\tau}=1)
\leq\sum_{i=1}^nP_0(E_i)=nP_0(E_1).
$$
For every integer $m\geq(1-\eta)\mun$, the event $C_1\leq\tau$ implies $|C_1-\Cd|\geq\Cd-\tau$. Lemma \ref{lem:null-conc} gives
$$
\begin{aligned}
&P_0(C_1\leq\tau\mid D_1=m) \leq2\exp\left(-\frac{(\Cd-\tau)^2m}{2}\right)\\
&\quad \leq2\exp\left(-\frac{(\Cd-\tau)^2(1-\eta)\mun}{2}\right).
\end{aligned}
$$
Averaging this conditional bound over $m$ proves the type-I error bound.

For the type-II error, fix $B\in\calB_{n,k_n}$ and choose $b\in B$. Since $E_b$ implies rejection,
$$
\begin{aligned}
&P_B(\psi_{n,\eta,\tau}=0) \leq P_B(E_b^c) \leq P_B(D_b<(1-\eta)\mun)\\
&\quad +P_B(D_b\geq(1-\eta)\mun,\ C_b>\tau).
\end{aligned}
$$
The degree of a bot is $\Bin(n-1,p_n)$, hence its mean is $\mun$. The binomial lower tail inequality gives
$$
P_B(D_b<(1-\eta)\mun)
\leq\exp\left(-\frac{\eta^2\mun}{2}\right).
$$
Conditional on $D_b=m\geq(1-\eta)\mun$, Lemma \ref{lem:bot-conc} gives
$$
\begin{aligned}
&P_B(C_b>\tau\mid D_b=m)\\
&\quad \leq P_B(|C_b-p_n|\geq\tau-p_n\mid D_b=m)\\
&\quad \leq2\exp\left(-\frac{(\tau-p_n)^2m}{4}\right),
\end{aligned}
$$
which is at most the second exponential term in the theorem. Average over $m$ and then maximize over $B$.
\end{proof}

The finite-sample bound yields the following asymptotic result.

\begin{corollary}\label{cor:power}
Assume $d\geq2$ is fixed, $p_n\to0$, $k_n$ is an integer satisfying $1\leq k_n\leq n$, and
$$
\frac{\mun}{\log n}\longrightarrow\infty.
$$
For any fixed $\eta\in(0,1)$ and any fixed $\tau\in(0,\Cd)$, the test $\psi_{n,\eta,\tau}$ is asymptotically powerful:
$$
P_0(\psi_{n,\eta,\tau}=1)
+
\max_{B\in\calB_{n,k_n}}P_B(\psi_{n,\eta,\tau}=0)
\longrightarrow0.
$$
In particular, the choice $\eta=1/2$ and $\tau=\Cd/2$ is valid for all sufficiently large $n$.
\end{corollary}

\begin{proof}
Since $p_n\to0$, eventually $r_n<\frac{1}{3}$ and $p_n < \tau$. In Theorem \ref{thm:nonasymp}, the coefficients multiplying $\mun$ in all three exponents are eventually bounded below by positive constants. The factor $n$ in the first term is overcome because $\mun$ dominates $\log n$; the other two terms have no prefactor depending on $n$. Thus all three vanish.
\end{proof}

The same thresholds also give an estimator of the bot set:
$$
\Bh_n=
\left\{i\in V_n:
D_i\geq\frac{\mun}{2}\;\text{and}\;C_i\leq\frac{\Cd}{2}\right\}.
$$
It is computed vertex by vertex and requires no global optimization.

\begin{theorem}[Exact localization by local clustering]\label{thm:localization}
Let $(k_n)$ be a deterministic integer sequence with $1\leq k_n=o(n)$. Assume $d\geq2$ is fixed, $p_n\to0$, and
$$
\frac{\mun}{\log n}\longrightarrow\infty.
$$
Then
$$
P_0(\Bh_n=\varnothing)\longrightarrow1.
$$
Moreover, the following convergence is uniform over all label sets of size $k_n$:
$$
\min_{B\in\calB_{n,k_n}}P_B(\Bh_n=B)
\longrightarrow1.
$$
\end{theorem}

\begin{proof}
Eventually $r_n < \frac{1}{3}$, $p_n < \Cd/2$, and $\mun/2\geq2$. The null claim is Theorem \ref{thm:nonasymp} with $\eta = 1/2$ and $\tau = \Cd/2$.

Fix $B\in\calB_{n,k_n}$. For each $b\in B$, the degree $D_b$ has distribution $\Bin(n-1,p_n)$, and hence
$$
P_B(D_b<\mun/2)\leq\exp(-\mun/8).
$$
Since $p_n\to0$, for all large $n$ one has $p_n\leq\Cd/4$. By Lemma \ref{lem:bot-conc}, for all $m\geq\mun/2$,
$$
\begin{aligned}
&P_B(C_b>\Cd/2\mid D_b=m) \leq2\exp\left(-\frac{(\Cd/4)^2m}{4}\right)\\
&\quad \leq2\exp\left(-\frac{\Cd^2\mun}{128}\right).
\end{aligned}
$$
A union bound over the $k_n\leq n$ bot vertices gives
$$
P_B(B\subset\Bh_n)
\geq1-n\exp(-\mun/8)-2n\exp\left(-\frac{\Cd^2\mun}{128}\right),
$$
and the right-hand side tends to one.

Now fix a genuine vertex $i \notin B$. Write
$$
M_i=|N(i)\cap (V_n\setminus B)|,
\qquad
H_i=|N(i)\cap B|.
$$
Thus $D_i=M_i+H_i$. Let $T_i^g$ be the number of edges induced by $N(i)\cap(V_n\setminus B)$, and let
$$
C_i^g=
\begin{cases}
T_i^g\binom{M_i}{2}^{-1},&M_i\geq2,\\
0,&M_i<2.
\end{cases}
$$
Fix $x \in \TT^d$ and a set $S\subset V_n\setminus(B\cup\{i\})$ with $|S| = m$. Given $X_i = x$ and $N(i)\cap(V_n\setminus B) = S$, the points $(X_j)_{j \in S}$ are independent and uniform on $B_{\TT}(x,r_n)$. Bot-incident edges are independent of these points and do not enter $C_i^g$. Thus, for $m \geq 2$, the conditional mean of $C_i^g$ is $\Cd$. The bounded-differences argument from Lemma \ref{lem:null-conc} gives

$$
\begin{aligned}
&P_B\left(C_i^g \leq \frac{3\Cd}{4} \mid X_i = x,\,N(i)\cap(V_n\setminus B) = S\right)\\
&\quad \leq 2\exp\left(-\frac{\Cd^2m}{32}\right).
\end{aligned}
$$

Averaging over $x$ and $S$ conditional on $M_i = m$ gives, for every feasible $2\leq m\leq n-k_n-1$,
$$
P_B\left( C_i^g \leq \frac{3\Cd}{4} \mid M_i = m\right) \leq 2 \exp\left(-\frac{\Cd^2m}{32}\right).
$$
Choose $\delta = 1/40$ and put
$$
\rho=\frac{\delta}{1/2-\delta}=\frac{1}{19}.
$$
If $D_i \geq \mun/2$, $H_i \leq \delta \mun$, and $C_i^g > 3\Cd/4$, then
$$
M_i = D_i - H_i \geq \left(\frac{1}{2}-\delta\right) \mun.
$$
Moreover,
$$
H_i \leq \rho M_i.
$$
For all sufficiently large $n$, the lower bound $M_i\geq(1/2-\delta)\mun$ also gives
$$
M_i \geq \frac{2(1+\rho)}{\rho}.
$$

The edges among the genuine neighbors already contribute $T_i^g$ triangles at $i$. Therefore
$$
C_i \geq C_i^g\frac{\binom{M_i}{2}}{\binom{M_i+H_i}{2}}.
$$
Since $H_i\leq\rho M_i$,
$$
\begin{aligned}
&\frac{\binom{M_i}{2}}{\binom{M_i+H_i}{2}} = \frac{M_i(M_i-1)}{(M_i+H_i)(M_i+H_i-1)}\\
&\quad \geq (1-M_i^{-1})(1+\rho)^{-2} \geq (1-\rho)(1+\rho)^{-2}.
\end{aligned}
$$
With $\rho = 1/19$ and $C_i^g > 3\Cd/4$,
$$
\begin{aligned}
&C_i >
\frac{3\Cd}{4}
\left(1-\frac{1}{19}\right)
\left(1+\frac{1}{19}\right)^{-2}\\
&\quad = \frac{513}{800}\Cd > \frac{3\Cd}{5} > \frac{\Cd}{2}.
\end{aligned}
$$
Therefore, a genuine vertex can pass the two thresholds defining $\Bh_n$ only if either $H_i>\delta\mun$, or $M_i\geq(1/2-\delta)\mun$ and $C_i^g\leq3\Cd/4$.

It remains to bound these two events. Here $H_i$ is $\Bin(k_n,p_n)$ with mean $\lambda_n=k_np_n$, and
$$
\frac{\lambda_n}{\mun}
=\frac{k_n}{n-1}=o(1).
$$
Eventually $\lambda_n\leq \delta\mun e^{-2}$, so the binomial upper-tail bound gives
$$
P_B(H_i>\delta\mun)
\leq\left(\frac{e\lambda_n}{\delta\mun}\right)^{\delta\mun}
\leq e^{-\delta\mun}.
$$
Hence
$$
P_B\left(\exists i\notin B:H_i>\delta\mun\right)
\leq n e^{-\delta\mun}\to0.
$$
For the second event, the conditional clustering bound gives
$$
\begin{aligned}
&P_B\left(M_i\geq\left(\frac{1}{2}-\delta\right)\mun,
\ C_i^g\leq\frac{3\Cd}{4}\right)\\
&\quad \leq
2\exp\left(-\frac{\Cd^2(1/2-\delta)\mun}{32}\right).
\end{aligned}
$$
A union bound over the genuine vertices shows that the probability of either exceptional event tends to zero. Hence
$$
P_B(\Bh_n\subset B)\to1
$$
uniformly over $B$. Together with $P_B(B\subset\Bh_n)\to1$, this gives
$$
\min_{B\in\calB_{n,k_n}}P_B(\Bh_n=B)\to1.
$$
\end{proof}

The proof gives a fixed gap between bot and genuine clustering values before thresholding.

\begin{corollary}[Uniform clustering margin]\label{cor:margin}
Under the assumptions of Theorem \ref{thm:localization}, fix $\gamma\in(0,\Cd/4)$. Then, uniformly over $B\in\calB_{n,k_n}$, the following statements hold with $P_B$-probability tending to one:
$$
D_b\geq\frac{\mun}{2}
\quad\text{and}\quad
C_b\leq\gamma
\qquad\text{for every }b\in B,
$$
and
$$
C_i\geq\frac{3\Cd}{5}
\qquad\text{for every }i\notin B\text{ such that }D_i\geq\frac{\mun}{2}.
$$
\end{corollary}

\begin{proof}
For bots, eventually $p_n\leq\gamma/2$. Lemma \ref{lem:bot-conc} and the degree lower-tail bound used in Theorem \ref{thm:localization} give
$$
\begin{aligned}
&P_B\left(\exists b\in B:D_b<\mun/2\text{ or }C_b>\gamma\right) \leq ne^{-\mun/8}\\
&\quad + 2n\exp\left(-\frac{\gamma^2\mun}{32}\right),
\end{aligned}
$$
which tends to zero.
For genuine vertices, use the two exceptional events from the proof of Theorem \ref{thm:localization}:
$$
H_i > \delta\mun
$$
and
$$
M_i\geq\left(\frac{1}{2}-\delta\right)\mun, \; C_i^g\leq\frac{3\Cd}{4},
$$
Outside these events, the deterministic estimate gives
$$
C_i> \frac{513}{800}\Cd > \frac{3\Cd}{5}.
$$
The same union bounds from the proof of Theorem \ref{thm:localization} show that these events occur with probability tending to zero, uniformly over $B$.
\end{proof}

\paragraph{Estimating the geometric benchmark.}
The null and recovery statements remain valid when $\Cd$ is replaced by a graph-measurable estimator $\widehat C_n$ satisfying
$$
\widehat C_n/\Cd\longrightarrow1
$$
in $P_0$-probability and uniformly in $P_B$-probability over $B\in\calB_{n,k_n}$. Fix $0<\varepsilon<1/5$. Whenever $|\widehat C_n/\Cd-1|\leq\varepsilon$, the threshold $\widehat C_n/2$ lies strictly between $2\Cd/5$ and $3\Cd/5$. For example, taking $\gamma=\Cd/5$ in Corollary \ref{cor:margin} separates all bots from the genuine candidates under $P_B$. Under $P_0$, false discoveries are controlled by Theorem \ref{thm:nonasymp} with $\tau=3\Cd/5$. This argument does not require an independent calibration sample. Bet et al. give a global clustering estimator that is consistent under both hypotheses \cite{bet2020detecting}, and adjacency-based dimension estimation has been studied in \cite{atamanchuk2024dimension}. We still assume that $\mun$ is known; joint finite-sample calibration of both thresholds is not addressed here.

A single vertex degree cannot provide this separation.

\begin{proposition}[Degree law at one vertex]\label{prop:degree-invisibility}
Fix $B\in\calB_{n,k_n}$. Under $P_B$, every fixed vertex $i\in V_n$, whether bot or genuine, satisfies
$$
D_i\sim\Bin(n-1,p_n).
$$
Under $P_0$, for every fixed vertex $i$,
$$
D_i\sim\Bin(n-1,p_n).
$$
Equivalently, the total variation distance between the degree laws of any prescribed vertex under $P_B$ and $P_0$ is zero. Hence any test based only on that one degree has the same distribution under both hypotheses.
\end{proposition}

\begin{proof}
If $i\in B$, its incident edges are independent Bernoulli variables with parameter $p_n$ by definition. If $i\notin B$, condition on $X_i=x$. Its edges to genuine vertices are independent Bernoulli variables with parameter $p_n$, since the other locations are independent and uniform. Its edges to bots have the same law and are independent of all locations. Hence in either case the $n-1$ incident indicators are independent Bernoulli variables with parameter $p_n$. The same argument applies under $P_0$. The degree laws and testing statement follow.
\end{proof}

This assertion concerns one degree at a time. The joint degree vector need not have the same law under the two hypotheses, and the proposition does not rule out tests based on the full degree sequence. It only explains why we condition on the degree of the vertex whose neighborhood is being examined.

The same concentration argument also works with deterministic degree and clustering thresholds. Let $a_n\geq2$ and $\tau_n$ be deterministic. Assume that there are constants $c_a>0$, $\eta_a\in(0,1)$, $\Delta>0$, and $c_0>0$ such that, for all sufficiently large $n$,

$$
\begin{gathered}
c_a\mun\leq a_n\leq(1-\eta_a)\mun,\\
\tau_n\leq\Cd-\Delta,\\
\tau_n-p_n\geq c_0.
\end{gathered}
$$
Define
$$
\psi_n^{a,\tau} = \one\{\exists i\in V_n:D_i\geq a_n,\ C_i\leq\tau_n\}.
$$

\begin{corollary}[Robust deterministic thresholds]\label{cor:det-threshold}
If $(k_n)$ is an integer sequence satisfying $1\leq k_n\leq n$ and $\mun/\log n\to\infty$, then $\psi_n^{a,\tau}$ is asymptotically powerful. More precisely, for all sufficiently large $n$,
$$
\begin{aligned}
&P_0(\psi_n^{a,\tau}=1) \leq 2n\exp\left(-\frac{(\Cd-\tau_n)^2a_n}{2}\right)\\
&\quad \leq 2n\exp\left(-\frac{\Delta^2c_a\mun}{2}\right),
\end{aligned}
$$
and
$$
\begin{aligned}
&\max_{B\in\calB_{n,k_n}}P_B(\psi_n^{a,\tau}=0) \leq \exp\left(-\frac{\eta_a^2\mun}{2}\right)\\
&\quad + 2\exp\left(-\frac{c_0^2a_n}{4}\right).
\end{aligned}
$$
Consequently,
$$
P_0(\psi_n^{a,\tau}=1) + \max_{B\in\calB_{n,k_n}}P_B(\psi_n^{a,\tau}=0) \to 0.
$$
\end{corollary}

\begin{proof}
The type-I bound follows by repeating the null argument from the proof of Theorem \ref{thm:nonasymp} with thresholds $a_n$ and $\tau_n$.

For the type-II bound, fix $B\in\calB_{n,k_n}$ and choose any $b \in B$. Since the event
$$
\{D_b \geq a_n,\ C_b \leq \tau_n\}
$$
implies $\psi_n^{a,\tau} = 1$, we have
$$
P_B(\psi_n^{a,\tau}=0) \leq P_B(D_b<a_n) + P_B(D_b\geq a_n,\ C_b>\tau_n).
$$
Because $D_b\sim\Bin(n-1,p_n)$ and $a_n\leq(1-\eta_a)\mun$, the binomial lower tail bound gives
$$
P_B(D_b<a_n) \leq \exp\left(-\frac{\eta_a^2\mun}{2}\right).
$$
For the clustering term, condition on $D_b=m\geq a_n$ and use Lemma \ref{lem:bot-conc}:
$$
\begin{aligned}
&P_B(C_b>\tau_n\mid D_b=m) \leq 2\exp\left(-\frac{(\tau_n-p_n)^2m}{4}\right)\\
&\quad \leq 2\exp\left(-\frac{c_0^2a_n}{4}\right).
\end{aligned}
$$
Average over $m$. Since $a_n\geq c_a\mun$ and $\mun/\log n\to\infty$, every error term tends to zero uniformly over $B$.
\end{proof}

The next graph-theoretic bound controls the effect of bot neighbors on the clustering of a genuine vertex.

\begin{lemma}[Deterministic neighborhood contamination inequality]\label{lem:det-contam}
Let $G$ be any finite simple graph, let $i$ be a vertex, and let the neighborhood of $i$ be decomposed as a disjoint union $N(i)=A\cup H$. Put $m=|A|$, $h=|H|$, and assume $m\geq2$. Let $T_A$ be the number of edges induced by $A$, and define $C_A=T_A\binom{m}{2}^{-1}$. Then
$$
C_i\geq C_A\frac{\binom{m}{2}}{\binom{m+h}{2}}.
$$
If $\rho\in(0,1)$, $h\leq\rho m$, and $m\geq 2(1+\rho)/\rho$, then
$$
C_i\geq C_A(1+\rho)^{-2}(1-\rho).
$$
In particular, along any sequence of graphs for which $m\to\infty$ and $h/m\to0$, one has $C_i \geq C_A (1 - o(1))$.
\end{lemma}

\begin{proof}
The full number of edges induced by $N(i)$ is at least the number of edges induced by $A$. Hence
$$
T_i\geq T_A=C_A\binom{m}{2}.
$$
Dividing by $\binom{m+h}{2}$ gives the first inequality. If $h\leq\rho m$, then
$$
\begin{aligned}
&\frac{\binom{m}{2}}{\binom{m+h}{2}} =\frac{m(m-1)}{(m+h)(m+h-1)} \geq\frac{m(m-1)}{m^2(1+\rho)^2}\\
&\quad =(1-m^{-1})(1+\rho)^{-2}.
\end{aligned}
$$
The lower bound on $m$ implies $1-m^{-1}\geq1-\rho$, proving the second inequality. For the asymptotic statement, use the first inequality and the exact ratio
$$
\frac{\binom{m}{2}}{\binom{m+h}{2}}
=\frac{1-m^{-1}}{(1+h/m)(1+(h-1)/m)}
\longrightarrow1
$$
whenever $m\to\infty$ and $h/m\to0$.
\end{proof}

Thus, if $C_A$ stays a fixed distance above the test threshold, adding a vanishing fraction of arbitrary neighbors cannot push $C_i$ below that threshold. In Theorem \ref{thm:localization}, $A$ is the set of genuine neighbors and $H$ is the set of bot neighbors.

The null concentration bound also yields vertex-level scores that account for the scan over all vertices. For $D_i\geq2$, define
$$
\begin{gathered}
\Pi_i=
1\wedge 2\exp\left(-\frac{D_i(\Cd-C_i)_+^2}{2}\right),\\
(x)_+=\max\{x,0\}.
\end{gathered}
$$
For $D_i<2$, set $\Pi_i=1$. A small $\Pi_i$ indicates unusually low clustering relative to $\Cd$.

\begin{proposition}[Finite-sample family-wise calibration]\label{prop:pvalues}
Suppose $r_n < \frac{1}{3}$. Under $H_0$, for every $\alpha \in (0,1)$,
$$
P_0\left(\min_{1\leq i\leq n}\Pi_i\leq\frac{\alpha}{n}\right)
\leq\alpha.
$$
Under $P_B$, let $b\in B$ be any bot. If $\gamma\in(0,\Cd)$ and $p_n<\gamma$, then
$$
\begin{aligned}
&P_B\left(\Pi_b>
2\exp\left(-\frac{\mun(\Cd-\gamma)^2}{4}\right)\right) \leq
\exp(-\mun/8)\\
&\quad +2\exp\left(-\frac{(\gamma-p_n)^2\mun}{8}\right).
\end{aligned}
$$
Thus a bot with degree at least $\mun/2$ and local clustering at most $\gamma$ has a score exponentially small in $\mun$.
\end{proposition}

\begin{proof}
Fix a vertex $i$ and a number $u\in(0,1)$. Conditional on $D_i=m\geq2$, the event $\Pi_i\leq u$ implies
$$
2\exp\left(-\frac{m(\Cd-C_i)_+^2}{2}\right)
\leq u.
$$
Therefore $C_i\leq\Cd$ and
$$
|C_i-\Cd|
\geq\left(\frac{2\log(2/u)}{m}\right)^{1/2}.
$$
Lemma \ref{lem:null-conc} gives
$$
P_0(\Pi_i\leq u\mid D_i=m)
\leq2\exp\left(-\frac{1}{2}m\frac{2\log(2/u)}{m}\right)=u.
$$
When $m<2$, one has $\Pi_i=1$, so $P_0(\Pi_i\leq u\mid D_i=m)=0$ because $u<1$. Averaging over $D_i$ gives $P_0(\Pi_i\leq u)\leq u$. Set $u=\alpha/n$ and apply the union bound.

For the bot bound, write
$$
q_n=2\exp\left(-\frac{\mun(\Cd-\gamma)^2}{4}\right).
$$
If $q_n\geq1$, then $\Pi_b\leq1\leq q_n$, so there is nothing to prove. Suppose $q_n<1$. On $\{D_b\geq\mun/2,\ C_b\leq\gamma\}$,
$$
\Pi_b
\leq2\exp\left(-\frac{D_b(\Cd-\gamma)^2}{2}\right)
\leq q_n.
$$
Hence $\Pi_b>q_n$ requires either $D_b<\mun/2$ or $D_b\geq\mun/2$ and $C_b>\gamma$. The first event has probability at most $e^{-\mun/8}$. For the second, degrees below two are irrelevant because $C_b=0<\gamma$. For every integer $m\geq2$ satisfying $m\geq\mun/2$, Lemma \ref{lem:bot-conc} gives
$$
\begin{aligned}
&P_B(C_b>\gamma\mid D_b=m) \leq2\exp\left(-\frac{(\gamma-p_n)^2m}{4}\right)\\
&\quad \leq2\exp\left(-\frac{(\gamma-p_n)^2\mun}{8}\right).
\end{aligned}
$$
Averaging over the degree completes the proof.
\end{proof}

The scores rank vertices without imposing a fixed degree cutoff. The asymptotic test in Theorem \ref{thm:nonasymp} keeps the cutoff because local clustering is unstable when a vertex has very few neighbors.

The same test can be expressed directly in terms of triangle counts. Define the signed wedge deficit
$$
W_i=\binom{D_i}{2}(\Cd-C_i)=\Cd\binom{D_i}{2}-T_i,
$$
with $W_i=0$ when $D_i<2$. Because $\Cd$ is real-valued, $W_i$ may be negative or non-integer. The following corollary rewrites the clustering condition in terms of $W_i$.

\begin{corollary}[Benchmark wedge deficit separation]\label{cor:wedge}
Let $\eta \in (0,1)$ and let $\tau$ satisfy $p_n<\tau < \Cd$. Assume $r_n < \frac{1}{3}$ and $(1 - \eta)\mun \geq 2$. For every vertex $i$ with $D_i\geq2$,
$$
C_i \leq \tau \; \Longleftrightarrow \; W_i \geq (\Cd - \tau) \binom{D_i}{2}.
$$
Hence,
$$
\begin{aligned}
&\{D_i\geq(1-\eta)\mun,\ C_i\leq\tau\}\\
&\quad = \left\{ D_i\geq(1-\eta)\mun, \ W_i\geq(\Cd-\tau)\binom{D_i}{2} \right\}.
\end{aligned}
$$
In particular, on this event,
$$
W_i \geq (\Cd - \tau)\binom{(1-\eta)\mun}{2},
$$
where $\binom{x}{2}$ denotes $x(x-1)/2$ when $x$ is not an integer.
Under $H_0$, the probability that the equivalent degree-normalized event in the event identity above holds for at least one vertex is at most the type-I bound in Theorem \ref{thm:nonasymp}. Under $P_B$, a fixed bot satisfies that same degree-normalized event with probability at least one minus the type-II bound in Theorem \ref{thm:nonasymp}.
\end{corollary}

\begin{proof}
For $D_i\geq2$,
$$
W_i = \binom{D_i}{2}(\Cd-C_i).
$$
Therefore $C_i \leq \tau$ is equivalent to
$$
W_i \geq (\Cd - \tau)\binom{D_i}{2}.
$$
Adding the condition $D_i\geq(1-\eta)\mun$ gives the displayed event identity. The lower bound with $\binom{(1-\eta)\mun}{2}$ follows from monotonicity of $x(x-1)/2$ for $x \geq 1$. The probability bounds now follow from Theorem \ref{thm:nonasymp}.
\end{proof}

The statistic uses adjacency information only. Build adjacency hash tables in $O(n+|E_n|)$ time, then test every pair in each neighborhood. This requires
$$
\sum_{i=1}^n\binom{D_i}{2}
=O\left(\sum_{i=1}^nD_i^2\right)
$$
adjacency queries. With standard hash tables, this is an expected-time bound. Under either hypothesis, Proposition \ref{prop:degree-invisibility} gives
$$
\Exp\left[\sum_{i=1}^nD_i^2\right]
=n\{\mun^2+\mun(1-p_n)\}.
$$
No latent positions or geometric reconstruction are needed. The theoretical scan uses $\Cd$ and $\mun$; $\Cd$ may be calibrated as described above.

\section{Discussion}
Local clustering separates bots from genuine vertices in this model even though a prescribed vertex has the same binomial degree law under both hypotheses. Conditional on any degree $m\geq2$, the mean local clustering coefficient is $\Cd$ under $P_0$ and $p_n$ at a bot under $P_B$. At a genuine vertex under $P_B$, the conditional mean is a mixture that converges to $\Cd$ when $k_n=o(n)$. Since $p_n\to0$ and $\Cd>0$, the gap between the bot and geometric centers remains of constant order.

We obtained an explicit non-asymptotic bound on the worst-case risk of the scan. If $d\geq2$ is fixed, $p_n\to0$, $\mun/\log n\to\infty$, and $(k_n)$ is any deterministic integer sequence with $1\leq k_n=o(n)$, the same local rule recovers the whole botnet exactly with probability tending to one, uniformly over its labels. The logarithmic factor pays for scanning $n$ neighborhoods after their clustering coefficients concentrate. This is a sufficient condition for the local statistic, not a minimax boundary for the full testing problem.

It remains to determine whether a matching lower bound holds for tests based only on local clustering. Finite-sample calibration is also open when both $\Cd$ and $\mun$ must be estimated, especially if the dimension grows with $n$; existing adjacency-based dimension estimators provide a starting point \cite{bet2020detecting,atamanchuk2024dimension}. Other extensions include soft distance kernels, repeated network observations, and adversarial bot attachments. A soft kernel replaces $\Cd$ by a kernel-dependent closure functional. Repeated observations add temporal information, while an adversarial extension must specify both a degree budget and the adversary's knowledge of the geometry.

Within the specified model, triangle closure can be used to test local consistency with latent geometry without observing the positions.

\section*{Data availability}
No data were created or analyzed in this study.

\appendix
\section{Auxiliary bounds and geometric facts}
We collect the standard estimates used above \cite{boucheron2013concentration,janson2000random}. Proofs are included to keep track of the constants.

\begin{lemma}[Chernoff lower tail]\label{app:chernoff-lower}
Let $X\sim\Bin(N,q)$ and $\mu=Nq$. For every $\eta\in(0,1)$,
$$
\Prob(X\leq(1-\eta)\mu)
\leq\exp\left(-\frac{\eta^2\mu}{2}\right).
$$
In particular,
$$
\Prob(X\leq\mu/2)\leq e^{-\mu/8}.
$$
\end{lemma}

\begin{proof}
For every $t>0$, Markov's inequality gives
$$
\begin{aligned}
&\Prob(X\leq(1-\eta)\mu) =\Prob(e^{-tX}\geq e^{-t(1-\eta)\mu})\\
&\quad \leq e^{t(1-\eta)\mu}\Exp e^{-tX}.
\end{aligned}
$$
Since $X$ is binomial,
$$
\begin{aligned}
&\Exp e^{-tX} =(1-q+qe^{-t})^N \leq\exp\{Nq(e^{-t}-1)\}\\
&\quad =\exp\{\mu(e^{-t}-1)\}.
\end{aligned}
$$
Choose $t=-\log(1-\eta)$. Then $e^{-t}=1-\eta$, and
$$
\Prob(X\leq(1-\eta)\mu)
\leq \exp\{\mu[(1-\eta)(-\log(1-\eta))-\eta]\}.
$$
The elementary inequality
$$
(1-\eta)(-\log(1-\eta))-\eta\leq-\eta^2/2
$$
for $\eta\in(0,1)$ proves the bound. Taking $\eta=1/2$ gives the last display.
\end{proof}

\begin{lemma}[Binomial upper tail]\label{app:chernoff-upper}
Let $X\sim\Bin(N,q)$ with mean $\lambda=Nq$. For every $x>\lambda$,
$$
\Prob(X\geq x)
\leq\left(\frac{e\lambda}{x}\right)^x.
$$
\end{lemma}

\begin{proof}
If $\lambda=0$, then $X=0$ almost surely, and the result is immediate for every $x>0$ with the convention $0^x=0$. Assume henceforth that $\lambda>0$.
For $t>0$, Markov's inequality gives
$$
\begin{aligned}
&\Prob(X\geq x) \leq e^{-tx}\Exp e^{tX} =e^{-tx}(1-q+qe^t)^N\\
&\quad \leq\exp\{-tx+\lambda(e^t-1)\}.
\end{aligned}
$$
Choose $t=\log(x/\lambda)>0$. Then
$$
\begin{aligned}
&\Prob(X\geq x) \leq\exp\{-x\log(x/\lambda)+x-\lambda\}\\
&\quad \leq\exp\{-x\log(x/(e\lambda))\} =\left(\frac{e\lambda}{x}\right)^x.
\end{aligned}
$$
\end{proof}

\begin{lemma}[McDiarmid's inequality]\label{app:mcdiarmid}
Let $Z_1,\ldots,Z_m$ be independent random variables taking values in arbitrary measurable spaces. Let $f$ be a measurable function such that replacing the $j$-th coordinate can change the value of $f$ by at most $c_j$. Then, for every $t>0$,
$$
\begin{aligned}
&\Prob\{|f(Z_1,\ldots,Z_m)-\Exp f(Z_1,\ldots,Z_m)|\geq t\}\\
&\quad \leq2\exp\left(-\frac{2t^2}{\sum_{j=1}^m c_j^2}\right).
\end{aligned}
$$
\end{lemma}

\begin{proof}
Let
$$
M_j=\Exp[f(Z_1,\ldots,Z_m)\mid Z_1,
\ldots,Z_j],
\qquad 0\leq j\leq m.
$$
Then $(M_j)_{j=0}^m$ is a martingale and
$$
f-\Exp f=\sum_{j=1}^m(M_j-M_{j-1}).
$$
Fix $j$ and condition on $Z_1,\ldots,Z_{j-1}$. The bounded difference assumption implies that, as $Z_j$ varies, the conditional range of $M_j$ is at most $c_j$. Hence, the conditional range of the martingale difference $M_j-M_{j-1}$ is at most $c_j$. Hoeffding's lemma gives, for every $s\in\RR$,

$$
\Exp\left[ \exp\left( s(M_j-M_{j-1})\right) \mid Z_1,\ldots,Z_{j-1} \right]
\leq \exp\left(\frac{s^2c_j^2}{8}\right).
$$
Optimizing the Chernoff bound over $s>0$ gives
$$
\Prob(f-\Exp f\geq t)
\leq\exp\left(-\frac{2t^2}{\sum_j c_j^2}\right).
$$
Apply the same argument to $-f$ and add the two tail bounds.
\end{proof}

\begin{lemma}[Torus lifting]\label{app:torus-lifting}
Let $0<r<1/3$ and let $x\in\TT^d$. If $y,z\in B_{\TT}(x,r)$, then there exist Euclidean representatives $u,v\in B(0,r)$ such that
$$
\distT(y,z) \leq r
\quad\Longleftrightarrow\quad
\|u-v\|_2\leq r.
$$
Consequently, for threshold events at radius $r$ involving only points in one torus ball of radius $r$, torus geometry agrees with Euclidean geometry after translation and lifting.
\end{lemma}

\begin{proof}
Translate so that $x=0$ and choose representatives $u,v\in B(0,r)$. Since $\|u-v\|_2\leq2r$, for every nonzero $z\in\ZZ^d$,
$$
\|u-v+z\|_2\geq\|z\|_2-\|u-v\|_2\geq1-2r>r.
$$
Thus no nonzero lattice shift can produce a distance at most $r$. The torus distance is at most $r$ exactly when the zero shift satisfies $\|u-v\|_2\leq r$.
\end{proof}

\begin{lemma}[Volume of the unit ball intersection]\label{app:lens}
Let $0\leq s\leq2$. The volume of the intersection of two unit balls in $\RR^d$ whose centers are at distance $s$ is
$$
\Vol\bigl(B(0,1)\cap B(se_1,1)\bigr)
=\thetad I_{1-s^2/4}\left(\frac{d+1}{2},\frac{1}{2}\right).
$$
For $0\leq s\leq1$, this is the expression used in Theorem \ref{thm:constant}.
\end{lemma}

\begin{proof}
By symmetry, the intersection is the union of two congruent spherical caps cut from the unit ball by the hyperplane perpendicular to the line of centers and located at distance $s/2$ from each center. The volume of one cap is
$$
\theta_{d-1}\int_{s/2}^{1}(1-t^2)^{(d-1)/2}\,dt.
$$
Therefore, the intersection volume equals
$$
2\theta_{d-1}\int_{s/2}^{1}(1-t^2)^{(d-1)/2}\,dt.
$$
Put $u=1-t^2$. Then $dt=-du/(2\sqrt{1-u})$, and the last display becomes
$$
\theta_{d-1}\int_0^{1-s^2/4}
 u^{(d-1)/2}(1-u)^{-1/2}\,du.
$$
Using
$$
\thetad=\frac{\pi^{d/2}}{\Gamma(d/2+1)},
\qquad
\theta_{d-1}=\frac{\pi^{(d-1)/2}}{\Gamma((d+1)/2)},
$$
and the beta identity
$$
B\left(\frac{d+1}{2},\frac{1}{2}\right)
=\frac{\Gamma((d+1)/2)\Gamma(1/2)}{\Gamma(d/2+1)},
$$
one obtains
$$
\theta_{d-1}B\left(\frac{d+1}{2},\frac{1}{2}\right)=\thetad.
$$
Dividing the incomplete beta integral by this beta normalizer gives the regularized incomplete beta function and proves the formula.
\end{proof}

\end{document}